\documentclass[11pt,a4paper]{amsart}
\usepackage{amsfonts,amscd,amssymb,amsmath,amsthm,hyperref,mathrsfs,xcolor,lscape,longtable}
\usepackage{microtype}
\newcommand{\ie}{i.e.\ }

\newcommand{\bb}[1]{\mathbb{#1}}

\newcommand{\C}{\bb{C}}
\newcommand{\Q}{\bb{Q}}

\newcommand{\Z}{\bb{Z}}

\newtheorem{thm}{Theorem}

\newtheorem{prop}{Proposition}

\theoremstyle{definition}
\newtheorem{rmk}{Remark}
\theoremstyle{definition}
\newtheorem{example}{Example}
\numberwithin{equation}{section}
\newcommand{\Sp}{\mathrm{Sp}}

\newcommand{\p}{\mathbb{P}}

\def\l{{\mathfrak{l}}}

\newcommand{\GL}{\mathrm{GL}}
\newcommand{\SL}{\mathrm{SL}}

\begin{document}
\date{\today}
\title[Symplectic Hypergeometric Groups]{Symplectic Hypergeometric Groups of degree six}
\author{Jitendra Bajpai, Daniele Dona, Sandip Singh and Shashank Vikram Singh}
\address{Mathematisches Institut, Georg-August-Universit\"at G\"ottingen, Germany}
\email{jitendra@math.uni-goettingen.de}
\address{Mathematisches Institut, Georg-August-Universit\"at G\"ottingen, Germany}
\email{daniele.dona@mathematik.uni-goettingen.de}
\address{Department of Mathematics, Indian Institute of Technology Bombay, Mumbai, India}
\email{sandip@math.iitb.ac.in}
\address{Department of Mathematics, Indian Institute of Technology Bombay, Mumbai, India}
\email{shashank@math.iitb.ac.in}
\begin{abstract}
Our computations show that there is a total of $40$ pairs of degree six coprime polynomials $f,g$ where $f(x)=(x-1)^6$, $g$ is a product of cyclotomic polynomials, $g(0)=1$ and $f,g$ form a primitive pair. The aim of this article is to determine whether the corresponding $40$ symplectic hypergeometric groups with a maximally unipotent monodromy follow the same dichotomy between arithmeticity and thinness that holds for  the $14$ symplectic hypergeometric groups corresponding to the pairs of degree four polynomials $f,g$ where $f(x)=(x-1)^4$ and $g$ is as described above. As a result we prove that at least $18$ of these $40$ groups are arithmetic in $\mathrm{Sp}(6)$.

In addition, we extend our search to all degree six symplectic hypergeometric groups.  We find that there is a total of $458$  pairs of polynomials (up to scalar shifts) corresponding to such groups. For $211$ of them, the absolute values of the leading coefficients of the difference polynomials $f-g$ are at most $2$ and the arithmeticity of the corresponding groups follows from Singh and Venkataramana, while the arithmeticity of one more hypergeometric group follows from Detinko, Flannery and Hulpke.

In this article, we show the arithmeticity of $160$ of the remaining $246$ hypergeometric groups.
\end{abstract}

\maketitle

\section{Introduction}\label{sec:intro}
A hypergeometric differential equation of order $n$ is an ordinary differential equation of order $n$ with three regular singular points. It is defined on the thrice punctured Riemann Sphere $\p^{1}(\C)\backslash \{ 0,1,\infty\}$. Let $\theta = z \frac{d}{dz}$ and $$\alpha = (\alpha_1, \ldots , \alpha_n) , \beta = ( \beta_1, \ldots , \beta_n)  \in \C^{n}.$$ We define the hypergeometric differential equation of order $n$ by 
\begin{equation}\label{hde}
[z(\theta + \alpha_1) \cdots (\theta + \alpha_n) - (\theta+\beta_1 -1)\cdots (\theta+\beta_{n -1} -1 )] u(z) =0\,.
\end{equation}

This has $n$ linearly independent solutions which can be explicitly  expressed as hypergeometric  functions of type ${}_{n} F_{n-1}$ around any point $z \in \p^{1}(\C)\backslash \{ 0,1,\infty\}$. For $\alpha=(\alpha_1, \ldots , \alpha_n) $ and $\beta=(\beta_1, \ldots , \beta_{n-1})$, we define 

$${}_{n} F_{n-1}(\alpha_1,\ldots, \alpha_n; \beta_1, \ldots, \beta_{n-1} | z)= \sum_{k=0}^{\infty} \frac{(\alpha_1)_{k} \ldots (\alpha_n)_{k} }{(\beta_1)_{k} \ldots (\beta_{n-1})_{k}} \frac{z^{k}}{k!}\,,$$ where $(\alpha)_{k}= \frac{\Gamma(\alpha+k)}{\Gamma(\alpha)}$. Then $n$-linearly independent solutions $u(z)$ of the equation (\ref{hde}) are defined by the functions
\begin{equation}
z^{1-\beta_j} {}_{n} F_{n-1}(1+\alpha_1-\beta_j,\ldots, 1+\alpha_n-\beta_j; 1+\beta_1-\beta_j, \ldots, \overline{1+\beta_j -\beta_j}, \ldots,  1+\beta_{n} -\beta_j | z)
\end{equation}
where $\overline{1+\beta_j -\beta_j}$ represents the omission of the term in the above expression.

Now it follows that the fundamental group $\pi_1$ of $\p^{1}(\C)\backslash \{ 0,1,\infty\}$ acts on the (local) solution space of the hypergeometric equation (\ref{hde}) and we get the monodromy representation $\rho:\pi_1\longrightarrow \GL(V)$ where $V$ is the $n$ dimensional solution space of the differential equation (\ref{hde}) on a small neighbourhood of a point $z_0$ (say) in $\p^{1}(\C)\backslash \{ 0,1,\infty\}$. The subgroup $\rho(\pi_1)$ of $\GL(V)$ is said to be the monodromy group of the hypergeometric differential equation (\ref{hde}). We also call it the hypergeometric group associated to the parameters $\alpha = (\alpha_1, \ldots , \alpha_n) , \beta = ( \beta_1, \ldots , \beta_n)  \in \C^{n}.$ 

Levelt \cite[Theorem 3.5]{BH} showed that if $\alpha_j-\beta_k\notin\Z$ for all $1\le j,k\le n$, then there exists a basis of the solution space of the hypergeometric equation with respect to which the hypergeometric group corresponding to the parameters $\alpha = (\alpha_1, \ldots, \alpha_n), \beta = (\beta_1, \ldots, \beta_n)  \in \C^{n}$ is the subgroup of $\GL_n(\C)$ generated by the companion matrices of the polynomials
\[f(x)=\prod_{j=1}^n(x-e^{2\pi i\alpha_j}),\quad g(x)=\prod_{j=1}^n(x-e^{2\pi i\beta_j})\] and any other hypergeometric group having the same parameters is a conjugate of this one. Note that the condition $\alpha_j-\beta_k\notin\Z$ for all $1\le j,k\le n$ ensures that the polynomials $f$ and $g$ do not have any common root.

Now we consider the case $n=6$ and the pair of polynomials $f,g$ that are products of cyclotomic polynomials, do not have any common root, form a primitive pair (that is, there do not exist polynomials $f_1,g_1\in\Z[x]$ so that $f(x)=f_1(x^k), g(x)=g_1(x^k)$ for $k\ge 2$), and $f(0)=g(0)=1$. Then, it follows from Beukers and Heckman \cite[Theorem 6.5]{BH} that the corresponding hypergeometric group $\Gamma(f,g)$ preserves a non-degenerate symplectic form $\Omega$ on $\Q^6$ and $\Gamma(f,g)$ is Zariski dense inside the corresponding symplectic group $\Sp_\Omega$. So in our case $\Gamma(f,g)\subseteq\Sp_\Omega(\Z)$ and we determine the pairs $f,g$ corresponding to which $\Gamma(f,g)$ has finite index in $\Sp_\Omega(\Z)$; whenever this occurs we call $\Gamma(f,g)$ arithmetic in the corresponding symplectic group.

Note that we made our count of all such pairs $f,g$ up to scalar shifts. By this we mean that it is equivalent to study the hypergeometric groups $\Gamma(f,g)$  and $\Gamma(f', g')$ when $f'(x)=f(-x)$ and $g'(x)=g(-x)$. This equivalence can be explained following the Remark 1.2 of~\cite{S17} by making an appropriate transition of $4 \times 4$ matrices into $6 \times 6$ matrices. We explain this using the following example: consider the pairs of polynomials $f(x)=\Phi_{1}(x)^{6}$, $g(x)=\Phi_3(x)\Phi_{6}(x)^{2}$ associated to the pairs of parameters $\alpha=(0,0,0,0,0,0)$ and $\beta=(\frac{1}{3}, \frac{2}{3},\frac{1}{6}, \frac{5}{6},\frac{1}{6}, \frac{5}{6} )$ and the pairs $f(-x)=\Phi_{2}(x)^{6}$, $g(-x)=\Phi_6(x)\Phi_{3}(x)^{2}$ associated to the pairs of parameters $\alpha'=(\frac{1}{2},\frac{1}{2},\frac{1}{2},\frac{1}{2},\frac{1}{2},\frac{1}{2})$ and $\beta'=(\frac{1}{3}, \frac{2}{3},\frac{1}{3}, \frac{2}{3},\frac{1}{6}, \frac{5}{6} )$.  The pairs $\alpha', \beta'$ and the pairs $\alpha, \beta$ can be transformed into one another by simply adding $\frac{1}{2}$ in each of their entries.  Here $\Phi_{n}(x)$ denotes the $n^{th}$-cyclotomic polynomial. The pairs $\alpha, \beta$, or equivalently the pairs $f,g$, respecting all the conditions discussed above will be our {\emph{qualified pairs}} to be considered for the study of this article, and we find $458$ such pairs. 

The arithmeticity of the degree six symplectic hypergeometric groups has been also investigated by Detinko, Flannery and Hulpke~\cite{DFH} and they have found one arithmetic group associated to pairs of  polynomials $f=\Phi_{3}(x) \Phi_5(x)$ and $g=\Phi_{14}(x)$. This is listed in Table 2 of~\cite{DFH}. For the complete list of their investigation see~\cite{H}. Notice that all  mentioned arithmetic groups in their list, except the groups associated to the pair of polynomials $f=\Phi_{3}(x) \Phi_5(x)$, $g=\Phi_{14}(x)$ and the pair $f'=\Phi_6(x) \Phi_{10}(x)$, $g'=\Phi_7(x)$ (which is simply a scalar shift of the pair $f, g$), are arithmetic by the criterion of Singh and Venkataramana~\cite[Theorem 1.1]{SV}.  
 
The following proposition easily follows from Singh and Venkataramana,  see \cite[Remark 5.1]{SV}.

\begin{prop}\label{Proposition}
 Let $f,g$ be a pair of degree $n$ polynomials which are products of cyclotomic polynomials, do not have any common roots, form a primitive pair and have the constant terms equal to $1$ (these conditions ensure that $n$ must be even). Let the leading coefficient of the difference polynomial $f-g$ has the absolute value $\ge 3$. Let $e_1,e_2,\ldots,e_n$ be the standard basis vectors of $\Q^n$ over $\Q$ and $I$ be the $n\times n$ identity matrix. Let $A,B$ be the companion matrices of the polynomials $f,g$, respectively, and $v=(A^{-1}B-I)(e_n)$. 
 
 If there exists an element $\gamma\in\Gamma(f,g)$ such that the three vectors $v, \gamma (v), \gamma^{-1}(v)$ are {\it linearly independent} and the coefficient of $e_n$ in $\gamma(v)$ is either $\pm 2$ or $\pm 1$, then the corresponding hypergeometric group $\Gamma(f,g)$ is arithmetic in the corresponding symplectic group.
\end{prop}

The above proposition is proved just by replacing either $A^k$ or $B^k$, depending on whether $\{v,A^k(v), A^{-k}(v)\}$ or $\{v,B^k(v), B^{-k}(v)\}$ is linearly independent (cf. \cite[Lemma 4.2]{SV}),  by $\gamma$ in the proof of  \cite[Theorem 1.1]{SV}. For the sake of completeness we provide a proof of the above proposition using \cite[Theorem 1.2]{SV} in Section \ref{ProofoftheProposition}.

It follows then that to show the arithmeticity of a symplectic hypergeometric group we only need to find an element $\gamma\in\Gamma(f,g)$ that satisfies the hypotheses of Proposition \ref{Proposition}. To apply this criterion we look at the hypergeometric groups $\Gamma(f,g)$ (where $f,g$ are products of cyclotomic polynomials) inside $\Sp(6)$ and find that there are in total $458$ hypergeometric groups satisfying the conditions of Beukers and Heckman \cite{BH} so that they are Zariski dense inside the corresponding symplectic groups (cf. Tables A, B, C and D). Out of these $458$ groups, there are $211$ (cf. Table C) satisfying the criterion of Singh and Venkataramana \cite[Theorem 1.1]{SV} and their arithmeticity follows. There are $247$ remaining groups (cf. Tables A, B and D) which do not satisfy the criterion of Singh and Venkataramana \cite[Theorem 1.1]{SV} and out of them there are $161$ (cf. Table A and Table B) which satisfy the hypotheses of Proposition \ref{Proposition} and their arithmeticity follows.  

\begin{rmk}
Note that  the linear independence condition in  Proposition~\ref{Proposition} is not reduntant: it is not always true that  for a $\gamma\in\Gamma(f,g)$ for which the coefficient of $e_n$ in  $\gamma(v)$ has absolute value 1 or 2, the three vectors $v, \gamma (v), \gamma^{-1}(v)$ are linearly independent. We have the following two examples.

\begin{example} In case $n=4$, let $$\alpha=\left(\frac{1}{2},\frac{1}{2},\frac{1}{3},\frac{2}{3}\right),\ \beta=\left(\frac{1}{4},\frac{1}{4},\frac{3}{4},\frac{3}{4}\right).$$
In this case the corresponding polynomials are \[f(x)=(x+1)^2(x^2+x+1)=x^4+3x^3+4x^2+3x+1,\ g(x)=(x^2+1)^2=x^4+2x^2+1.\]
Now, if we denote, respectively, by $A$ and $B$ the companion matrices of $f$ and $g$, then $v=(3,2,3,0)$ and for $\gamma=BA$, the coefficient of $e_4$ in $\gamma(v)$ is $2$ but the vectors $v, \gamma (v), \gamma^{-1}(v)$ are {\it not} linearly independent.
\end{example}

\begin{example}
 In case $n=6$, let $$\alpha=\left(\frac{1}{2},\frac{1}{2},\frac{1}{2},\frac{1}{2},\frac{1}{6},\frac{5}{6}\right),\ \beta=\left(\frac{1}{9},\frac{2}{9},\frac{4}{9},\frac{5}{9},\frac{7}{9},\frac{8}{9}\right).$$
In this case the corresponding polynomials are \[f(x)=(x+1)^4(x^2-x+1)={x}^{6}+3\,{x}^{5}+3\,{x}^{4}+2\,{x}^{3}+3\,{x}^{2}+3\,x+1,\ g(x)=x^6+x^3+1.\]
Now, if we denote, respectively, by $A$ and $B$ the companion matrices of $f$ and $g$, then $v=(3,3,1,3,3,0)$ and for $\gamma=B^2A$, the coefficient of $e_6$ in $\gamma(v)$ is $1$ but the vectors $v, \gamma (v), \gamma^{-1}(v)$ are {\it not} linearly independent.
 
\end{example}

Therefore, we cannot drop the linear independence condition from the above proposition if we want to use the method of the proof of \cite[Theorem 1.1]{SV}.
 \end{rmk}
 
One of the starting motivations behind the work of this article lies in an attempt to answer a question asked by N. Katz during the workshop on ``Thin Groups and Super Approximation" held at IAS Princeton in March 2016, where the first and the third author were among the participants. He asked whether the degree six symplectic hypergeometric groups with a maximally unipotent monodromy follow the same pattern as the 14 degree four symplectic hypergeometric groups with a maximally unipotent monodromy: we know in fact from \cite{S15S, SV} that 7 of the 14 degree four groups are arithmetic, and the other 7 are thin by \cite{BT}, and one may wonder whether a similar dichotomy occurs in this particular family of degree six symplectic hypergeometric groups. 

We  summarize our effort to answer the question above in the following theorem.

\begin{thm}\label{thm1}
There are $40$ degree six symplectic  hypergeometric groups with a maximally unipotent monodromy, listed in Table A, out of which at least $18$ are arithmetic.
\end{thm}
  
In addition to these $18$ arithmetic hypergeometric groups, we extend our search to the remaining hypergeometric groups: one of them is known to be arithmetic by~\cite{DFH}, and we are able to find $142$ more. More precisely, we conclude the following.
 
 \begin{thm}\label{thm2}
 The $143$ hypergeometric groups appearing in Table B in Section~\ref{arithmeticexamples} are arithmetic.
 \end{thm}
 
The arithmeticity of the groups mentioned in Theorem~\ref{thm1} and Theorem~\ref{thm2}  above follows from Proposition \ref{Proposition}. 
 
\section{Proof of Proposition \ref{Proposition}}\label{ProofoftheProposition}
It follows that the hypergeometric group $\Gamma(f,g)$ preserves a non-degenerate symplectic form $\Omega$ and $\Gamma(f,g)\subseteq\Sp_\Omega(\Z)$ is a Zariski dense subgroup (cf. \cite[Theorem 6.5]{BH}). Hence to use \cite[Theorem 1.2]{SV} we need to find three transvections $C_1,C_2,C_3\in \Gamma(f,g)$ and vectors $w_1,w_2,w_3\in\Z^n$ so that the set $\{w_1,w_2,w_3\}$ is linearly independent, $\Z w_i=(C_i-1)(\Z^n),\ \forall 1\le i\le 3$ and $\Omega(w_i,w_j)\neq 0$ for some $1\le i,j\le 3$. With these conditions it follows that each $C_i$ maps $W$, the subspace spanned by the vectors $w_1,w_2,w_3$, into itself and then we also need to show  that the group generated by the restrictions $C_1|_W,C_2|_W,C_3|_W$ contains a nontrivial element of the unipotent radical of $\Sp_W$.

We consider the following transvections $C_1=C=A^{-1}B$, $C_2=\gamma^{-1}C\gamma$, $C_3=\gamma C\gamma^{-1}$ and $w_1=v=(C-1)(e_n)$, $w_2=\gamma ^{-1}(v)$ and $w_3=\gamma (v)$. Note that $(C_1-1)(w)=\lambda v$ for all $w\in\Z^n$ and for some $\lambda\in\Z$, and then it follows that for each $1\le i\le 3$, $(C_i-1)(w)=\lambda_jw_i$ for all $w\in\Z^n$ and for some $\lambda_i\in\Z$. Hence the condition that $(C_i-1)(\Z^n)=\Z w_i, \forall 1\le i\le 3$ is satisfied. Also, it is part of the hypotheses of the proposition that the vectors $w_1,w_2,w_3$ are linearly independent. 

Just by using the invariance of $\Omega$ under the action of $C$ and its non-degeneracy, we find that $\Omega(v,e_j)=0, \forall 1\le j\le n-1$ and $\Omega(v,e_n)\neq 0$. If $c$ is the coefficient of $e_n$ in $\gamma (v)$, $\Omega(v,\gamma^{-1} v)=\Omega(\gamma v,v)=-c\Omega(v,e_n)\neq 0$ implies that $\Omega(w_1,w_2)\neq0$.

Now, we consider the $3$ dimensional subspace $W=\sum_{j=1}^3\Q w_j$ and show that the group generated by the restrictions of the $C_i$ (for $1\le i\le 3$) to $W$ contains a nontrivial element of the unipotent radical of the symplectic group of $W$. Since $\dim W=3$ (an odd number), the restriction $\Omega|_W$ of $\Omega$ on $W$ is degenerate and $W\cap W^\perp$ is one dimensional as $\Omega(w_1,w_2)\neq 0$. Let $e\in W$ be a vector such that $W\cap W^\perp=\left<e\right>$. Note that $e$ cannot be written as linear combination of $w_1$ and $w_2$, and hence the set $\{e,w_1,w_2\}$ is linearly independent and gives a basis of $W$. With respect to this basis, $\Sp(W)=\Q^2\rtimes\SL_2(\Q)$ can be realized as \[\left\{\begin{pmatrix}1&u_1&u_2\\0&a_1&a_2\\0&b_1&b_2
\end{pmatrix}: u_1,u_2,a_1,a_2,b_1,b_2\in\Q, a_1b_2-a_2b_1=1                                                                                                                                                                                                                                                                                                                                                                                                                                                                                                                                                                                                                                                                                                                            \right\}\]and \[\left\{\begin{pmatrix}1&u_1&u_2\\0&1&0\\0&0&1
\end{pmatrix}: u_1,u_2\in\Q                                                                                                                                                                                                                                                                                                                                                                                                                                                                                                                                                                                                                                                                                                                            \right\}\]is the unipotent radical of $\Sp(W)$. 

Now, we only need to check that 
\[\begin{pmatrix}1&u_1&u_2\\0&1&0\\0&0&1
\end{pmatrix}\in\left<C_1|_W,C_2|_W,C_3|_W\right>\]for some $(u_1,u_2)\neq (0,0)\in\Z^2$.

Since $C_j$ is unipotent and $C_j(e)\in W\cap W^\perp$, it follows that $C_j(e)=e$, $\forall 1\le j\le 3$. 

By an easy check we find that $C_1(w_1)=w_1$, $C_1(w_2)=C(\gamma^{-1}(v))=\gamma^{-1}(v)-cv=w_2-cw_1$ and it follows that
\[C_1|_W=\begin{pmatrix}
          1&0&0\\0&1&-c\\0&0&1
         \end{pmatrix}.
\]

Also, $C_2(w_1)=\gamma^{-1}C\gamma(v)=\gamma^{-1}(\gamma(v)+cv)=v+c\gamma^{-1}(v)=w_1+cw_2$, $C_2(w_2)=\gamma^{-1}C\gamma(\gamma^{-1}(v))=\gamma^{-1}(v)=w_2$ and it follows that
\[C_2|_W=\begin{pmatrix}
          1&0&0\\0&1&0\\0&c&1
         \end{pmatrix}.
\]
Note that, for $c=\pm1,\pm2$, the two matrices \[
\begin{pmatrix}
          1&-c\\0&1
         \end{pmatrix}, \begin{pmatrix}
          1&0\\c&1
         \end{pmatrix}
\]
generate a finite index subgroup of $\SL_2(\Z)$. 

Now, we write the matrix representation of $C_3|_W$ with respect to the basis $\{e,w_1,w_2\}$. For, $C_3(w_1)=\gamma C\gamma^{-1}(v)=\gamma C(-ce_n+v')=\gamma(-c(e_n+v)+v')=\gamma((-ce_n+v')-cv)=\gamma(\gamma^{-1}(v)-cv)=v-c\gamma (v)=w_1-cw_3$ where $v'$ is a linear combination of the vectors $e_1,e_2,\ldots,e_{n-1}$ and hence it is fixed under the action of $C$. If we write $C_3(w_1)=l_1 e+l_2w_1+l_3w_2$ for some $l_1,l_2,l_3\in\Q$, then it follows that $l_1\neq0$ as $w_1,w_2,w_3$  are linearly independent. Now, $C_3(w_2)=m_1e+m_2w_1+m_3w_2$ for some $m_1,m_2,m_3\in\Q$. Then, it follows that

\[C_3|_W=\begin{pmatrix}
          1&l_1&m_1\\0&l_2&m_2\\0&l_3&m_3
         \end{pmatrix}
\]with $l_1\neq0$.

Since $C_3$ is unipotent and for $c=\pm1$ or $\pm2$, $\begin{pmatrix}
          1&-c\\0&1
         \end{pmatrix}, \begin{pmatrix}
          1&0\\c&1
         \end{pmatrix}$ generate a finite index subgroup of $\SL_2(\Z)$, the $2\times 2$ matrix $u=\begin{pmatrix}
          l_2&m_2\\l_3&m_3
         \end{pmatrix}$ is unipotent and hence there exists an integer $m$ such that
         \[u^m\in \left<\begin{pmatrix}
          1&-c\\0&1
         \end{pmatrix}, \begin{pmatrix}
          1&0\\c&1
         \end{pmatrix}\right>\]and there exists $h\in \langle C_1|_W,C_2|_W\rangle$ such that \[(C_3|_W)^m\cdot h^{-1}=\begin{pmatrix}
                                           1&t_1&t_2\\0&1&0\\0&0&1\end{pmatrix}
\in \langle C_1|_W,C_2|_W,C_3|_W\rangle \]where $(t_1,t_2)=(l_1,m_1)(1+u+u^2+\cdots+u^{m-1})\neq (0,0)$ (since $u$ is unipotent and hence $1+u+u^2+\cdots+u^{m-1}$ is nonsingular; and $l_1\neq0$). Thus the element $(C_3|_W)^m\cdot h^{-1}$ of the group $\left< C_1|_W,C_2|_W,C_3|_W\right>$ is a non-trivial element of the unipotent radical of $\Sp_W$ and it now follows from \cite[Theorem 1.2]{SV} that the group $\Gamma(f,g)$ satisfying the hypotheses of Proposition \ref{Proposition}, is arithmetic. \qed

\section{Sage Code}\label{se:sage}

 In this section we present the program that aided us in detecting arithmetic hypergeometric groups. The program is written in SageMath, version 8.9~\cite{sage}: the computations are quite elementary and could have been performed by programs in other languages as well but SageMath is an open-source which is why we chose to use it.

The program is designed to take two polynomials $f,g$ and an integer $k$, and find whether there exists some $\gamma\in\Gamma(f,g)$ satisfying the hypotheses of Proposition~\ref{Proposition} and that can be written as a product of at most $k$ matrices in $\{A,B,A^{-1},B^{-1}\}$, where $A$ and $B$ are the companion matrices of $f$ and $g$ respectively. Example values have already been inserted, with $f=\Phi_{1}(x)^{6}$, $g=\Phi_{3}(x)^{2}\Phi_{6}(x)$ and $k=9$, corresponding to the parameters appearing in entry $17$ of Table A; an interested reader only needs to modify these values to use the code for themselves (see the lines immediately after ``\verb|# Here the main program starts|''). The code is commented throughout, to improve legibility and verifications.

\begin{verbatim}
#####
# Given an integer k and two polynomials f,g of the same degree (say n),
# the program takes their companion matrices A,B and the vector
# v=(A^(-1)B-I)e_n, then it finds out whether there is
# a product M of at most k matrices in {A,B,A^(-1),B^(-1)}
# such that the n-th entry of the vector Mv is in {+1,-1,+2,-2}
# and the three vectors M^(-1)v,v,Mv are linearly independent.
#####
# The following subroutine converts the vector vf of the coefficients
# of a polynomial f into the companion matrix of f.
def companion_internal(vf):
    le=len(vf)-1
    M=matrix(le,le)
    M[0,le-1]=-vf[0]
    i=1
    while i<le:
        M[i,i-1]=1
        M[i,le-1]=-vf[i]
        i+=1
    return M
#####
# The following subroutine converts a polynomial func
# into its companion matrix.
def companion(func):
    return companion_internal(func.list())
#####
# The following subroutine takes two nxn-matrices A,B and returns
# the vector v=(A^(-1)B-I)e_n.
def othervec(A,B,n):
    en=vector([0]*n)
    en[n-1]=1
    return (A^(-1)*B-matrix.identity(n))*en
#####
# The following subroutine checks if the last entry of a vector v
# is in {+1,-1,+2,-2}.
def checklast(v):
    entry=v[len(v)-1]
    if entry==1 or entry==-1 or entry==2 or entry==-2:
        return True
    else:
        return False
#####
# The following subroutine checks if M^(-1)v,v,Mv are
# linearly independent (in Q). It returns True if they
# are, and False if they are not.
def independent(M,v):
    E=QQ^len(v)
    return not E.are_linearly_dependent([M^(-1)*v,v,M*v])
#####
# The following subroutine tries both the "checklast" subroutine
# on a vector v and the "independent" subroutine on v and a matrix M.
# If both are True, it returns [True,M]; if at least one is False,
# it returns [False].
def tryone(M,v,s):
    if checklast(M*v):
        if independent(M,v):
            return [True,M,s]
        return [False]
    return [False]
#####
# The following is the main subroutine: given an integer k, a vector v,
# 5 matrices A,B,C,D,M (with C=A^(-1),D=B^(-1)), a string s
# of A,B,A^(-1),B^(-1) corresponding to M, and an indicator that says
# what the last matrix in the decomposition of M was, first it checks
# whether M and v satisfy the two desired properties; then if they
# do not and |s|<k, it calls recursively the same subroutine for
# M*A,M*B,M*A^(-1),M*B^(-1) (actually, only 3 of them, the ones that
# do not involve a X*X^(-1) at the end). If at any point there is one M
# that satisfies the two desired properties, it returns [True,M,s];
# if there is not, it returns [False].
def tryall(s,k,A,B,C,D,M,v,lastguy):
    check=tryone(M,v,s)
    if len(s)<k:
        if lastguy==0:
            Next=[['A',A],['B',B],['A^(-1)',C],['B^(-1)',D]]
        elif lastguy=='A':
            Next=[['A',A],['B',B],['B^(-1)',D]]
        elif lastguy=='B':
            Next=[['A',A],['B',B],['A^(-1)',C]]
        elif lastguy=='A^(-1)':
            Next=[['B',B],['A^(-1)',C],['B^(-1)',D]]
        elif lastguy=='B^(-1)':
            Next=[['A',A],['A^(-1)',C],['B^(-1)',D]]
        else:
            return 'Error'
        indnext=0
        while indnext<len(Next):
            if check==[False]:
                check=tryall(s+[Next[indnext][0]],k,A,B,C,D, \
                M*Next[indnext][1],v,Next[indnext][0])
            indnext+=1
    return check
#####
# Here the main program starts.
# Define f here.
f=cyclotomic_polynomial(1)^6
# Define g here.
g=cyclotomic_polynomial(3)^2*cyclotomic_polynomial(6)
# Define k here: this is the maximal length that one wants to check.
k=9
# Printing the input.
print('== Input ==')
print('Polynomial f:')
print(f)
print('Polynomial g:')
print(g)
A=companion(f)
B=companion(g)
C=A^(-1)
D=B^(-1)
v=othervec(A,B,f.degree(x))
print('Matrix A:')
print(A)
print('Matrix B:')
print(B)
print('Matrix A^(-1):')
print(C)
print('Matrix B^(-1):')
print(D)
print('Vector v:')
print(v)
print('We try up to length',k,'and see if there is a matrix that works.')
# Here the main subroutine is called.
final=tryall([],k,A,B,C,D,matrix.identity(f.degree(x)),v,0)
# If "final" is [False], it means that there was no product of at most
# k instances of A,B,A^(-1),B^(-1) such that the conditions are satisfied.
# Otherwise, "final" is [True,M,s] where M is the product matrix itself
# and s is a string of A,B,A^(-1),B^(-1) that represents M.
print(' ')
print('== Output ==')
if final[0]==False:
    print('There is no word of at most',k,'letters that respects all conditions.')
elif final[0]==True:
    print('There is a word respecting all conditions!')
    s=final[2]
    le=len(s)
    strin=''
    i=0
    while i<le:
        strin=strin+s[i]
        i+=1
    print('Product:')
    print(strin)
    print('Matrix M to which it corresponds:')
    print(final[1])
    print('Vectors Mv,v,M^(-1)v:')
    print(final[1]*v)
    print(v)
    print(final[1]^(-1)*v)
else:
    print('There is some error in this code.')
\end{verbatim}

\section{Table A:  hypergeometric groups in $\Sp(6)$ with $\alpha=(0,0,0,0,0,0)$}\label{sec:mum}
In this section, we consider an interesting family  of hypergeometric groups  with  a maximally unipotent monodromy, \ie with respect to the pairs $f,g$ of degree six polynomials where $f=(x-1)^6$ (\ie $\alpha=(0,0,0,0,0,0)$ in this case), $g$ is a product of cyclotomic polynomials and satisfying $g(0)=1, g(1)\neq 0$ and $f,g$ form a primitive pair.  Our computations show that there are 40 such qualified pairs which we list in the table below.

In the table below the second column records all the possible $\beta=(\beta_1,\beta_2,\ldots,\beta_6)$ which determines $g$ as $g=\prod_{j=1}^6(x-e^{2\pi i\beta_j})$, the third column keeps track of the absolute value of the leading coefficient of the difference polynomial $f-g$, denoted by $|lc(f-g)|$, the fourth column describes the vector $v=(A^{-1}B-I)(e_6)$ with $e_6=(0,0,0,0,0,1)$ and in the fifth column we provide a $\gamma\in \Gamma(f,g)$ for which the hypotheses of Proposition \ref{Proposition} are satisfied (when we know such a $\gamma$).  In the last column the reader will find the $18$ groups whose arithmeticity has been proved in this article, are marked with ``Yes", and the $22$ examples  whose arithmeticity or thinness  is unknown, are marked with ``??".  Notice that for all the $40$ groups listed below $|lc(f-g)| \geq 3$, so the criterion of \cite[Theorem 1.1]{SV} cannot be applied in these cases. However,  the arithmeticity of the $18$ of the 40 hypergeometric groups  follows from Proposition~\ref{Proposition}.

It has been expected that many of the examples from this family will correspond to the mirrors of
true Calabi-Yau 5-folds, for details see~\cite{GMP}. One such interesting example will be the septic case, \ie the parameter $\beta=\big(\frac{1}{7},\frac{2}{7},\frac{3}{7},\frac{4}{7},\frac{5}{7},\frac{6}{7}\big)$. The reader may think of it as interesting as  the quintic case of degree four associated to the parameters $\alpha=(0,0,0,0)$ and $\beta=(\frac{1}{5},\frac{2}{5},\frac{3}{5},\frac{4}{5}\big)$. These are also known in the literature as the members of the Dwork family, for quick details  see Appendix-I and III in~\cite{Katz} where one can think  about the quintic and septic cases in particular while keeping the values of $d=n=5$ and $d=n=7$, respectively, in~\cite{Katz}. 
{\begin{center}
\scriptsize\renewcommand{\arraystretch}{2}
\begin{longtable}{|c|c|c|c|c|c|}
\hline
 S.No. & $\beta$ & $|lc(f-g)|$ & $v$  & $\gamma$ & Arithmetic \\                                                                                                                                                        
\hline
\hline
1  & $\big(\frac{1}{2},\frac{1}{2},\frac{1}{2},\frac{1}{2},\frac{1}{2},\frac{1}{2}\big)$   & 12 & $(-12, 0, -40, 0, -12, 0)$&  ?? & ??  \\
\hline
2 &   $\big(\frac{1}{2},\frac{1}{2},\frac{1}{2},\frac{1}{2},\frac{1}{3},\frac{2}{3}\big)$  & 11 & $(-11, 4, -34, 4, -11, 0)$& ?? & ??\\
\hline
3 &   $\big(\frac{1}{2},\frac{1}{2},\frac{1}{2},\frac{1}{2},\frac{1}{4},\frac{3}{4}\big)$ & 10 & $(-10, 8, -28, 8, -10, 0)$& ?? &  ??  \\
\hline
4 &   $\big(\frac{1}{2},\frac{1}{2},\frac{1}{2},\frac{1}{2},\frac{1}{6},\frac{5}{6}\big)$  & 9 & $(-9, 12, -22, 12, -9, 0)$ & ?? & ??  \\
\hline
5 &   $\big(\frac{1}{2},\frac{1}{2},\frac{1}{3},\frac{1}{3},\frac{2}{3},\frac{2}{3}\big)$ & 10 &$(-10, 7, -30, 7, -10, 0)$  & ?? & ?? \\
\hline
6 &  $\big(\frac{1}{2},\frac{1}{2},\frac{1}{3},\frac{2}{3},\frac{1}{4},\frac{3}{4}\big)$ & 9 & $(-9, 10, -26, 10, -9, 0)$& ?? &??  \\
\hline
7 &  $\big(\frac{1}{2},\frac{1}{2},\frac{1}{3},\frac{2}{3},\frac{1}{6},\frac{5}{6}\big)$ & 8 & $(-8, 13, -22, 13, -8, 0)$& ?? & ??\\
\hline
8 & $\big(\frac{1}{2},\frac{1}{2},\frac{1}{4},\frac{1}{4},\frac{3}{4},\frac{3}{4}\big)$ & 8 & $(-8, 12, -24, 12, -8, 0)$ & ?? &  ?? \\
\hline
9 & $\big(\frac{1}{2},\frac{1}{2},\frac{1}{4},\frac{3}{4},\frac{1}{6},\frac{5}{6}\big)$ & 7 & $(-7, 14, -22, 14, -7, 0)$& ?? & ?? \\  
\hline
10 & $\big(\frac{1}{2},\frac{1}{2},\frac{1}{5},\frac{2}{5},\frac{3}{5},\frac{4}{5}\big)$ & 9 & $(-9, 11, -24, 11, -9, 0)$& ?? & ??  \\ 
 \hline
11 &   $\big(\frac{1}{2},\frac{1}{2},\frac{1}{6},\frac{1}{6},\frac{5}{6},\frac{5}{6}\big)$ & 6 & $(-6, 15, -22, 15, -6, 0)$& ?? & ??  \\  
\hline
12 &    $\big(\frac{1}{2},\frac{1}{2},\frac{1}{8},\frac{3}{8},\frac{5}{8},\frac{7}{8}\big)$ &  8 & $(-8, 14, -20, 14, -8, 0)$ & ?? & ?? \\  
\hline
13 & $\big(\frac{1}{2},\frac{1}{2},\frac{1}{10},\frac{3}{10},\frac{7}{10},\frac{9}{10}\big)$  &  7 &$(-7, 15, -20, 15, -7, 0)$ & ?? & ??\\  
\hline
14 & $\big(\frac{1}{2},\frac{1}{2},\frac{1}{12},\frac{5}{12},\frac{7}{12},\frac{11}{12}\big)$ & 8 & $(-8, 15, -18, 15, -8, 0)$ & ?? & ?? \\ 
 \hline
15 & $\big(\frac{1}{3},\frac{1}{3},\frac{1}{3},\frac{2}{3},\frac{2}{3},\frac{2}{3}\big)$  & 9 & $(-9, 9, -27, 9, -9, 0)$ & ?? & ?? \\  
\hline
16 & $\big(\frac{1}{3},\frac{1}{3},\frac{2}{3},\frac{2}{3},\frac{1}{4},\frac{3}{4}\big)$ & 8 & $(-8, 11, -24, 11, -8, 0)$ & ?? & ??  \\ 
 \hline
17 &  $\big(\frac{1}{3},\frac{1}{3},\frac{2}{3},\frac{2}{3},\frac{1}{6},\frac{5}{6}\big)$  & 7 &$(-7, 13, -21, 13, -7, 0)$ &  $A^2 B A^{-1}B^{4}A$ & Yes\\ 
 \hline
18 &  $\big(\frac{1}{3},\frac{2}{3},\frac{1}{4},\frac{1}{4},\frac{3}{4},\frac{3}{4}\big)$ & 7 & $(-7, 12, -22, 12, -7, 0)$ &$AB^{-1}A^3 B^3 A B^{-3}$ &  Yes\\  
\hline
19 & $\big(\frac{1}{3},\frac{2}{3},\frac{1}{4},\frac{3}{4},\frac{1}{6},\frac{5}{6}\big)$ & 6 &$(-6, 13, -20, 13, -6, 0)$ & $AB^2 AB^5 A B^{-1}A B^{-2}$ &  Yes  \\  
\hline
20 &  $\big(\frac{1}{3},\frac{2}{3},\frac{1}{6},\frac{5}{6},\frac{1}{6},\frac{5}{6}\big)$ & 5 & $(-5,13,-19,13,-5,0)$& $B^{3}$& Yes  \\ 
 \hline
21 &  $\big(\frac{1}{3},\frac{2}{3},\frac{1}{5},\frac{2}{5},\frac{3}{5},\frac{4}{5}\big)$  & 8 & $(-8, 12, -23, 12, -8, 0)$&?? & ?? \\  
\hline
22 & $\big(\frac{1}{3},\frac{2}{3},\frac{1}{8},\frac{3}{8},\frac{5}{8},\frac{7}{8}\big)$ & 7 & $(-7,14,-20,14,-7,0)$& $AB^{5}$ & Yes\\  
\hline
23 & $\big(\frac{1}{3},\frac{2}{3},\frac{1}{10},\frac{3}{10},\frac{7}{10},\frac{9}{10}\big)$  & 6 & $(-6, 14, -19, 14, -6, 0)$& $B^{6}A^{2}B^{4}A^{-1}$ &  Yes\\ 
 \hline
24 &  $\big(\frac{1}{3},\frac{2}{3},\frac{1}{12},\frac{5}{12},\frac{7}{12},\frac{11}{12}\big)$ &  7 &$(-7, 15, -19, 15, -7, 0)$ &?? & ?? \\  
\hline
25 & $\big(\frac{1}{4},\frac{1}{4},\frac{1}{4},\frac{3}{4},\frac{3}{4},\frac{3}{4}\big)$ & 6 &$(-6,12,-20,12,-6,0)$ & $B^{2}A$ & Yes  \\ 
 \hline
26 &  $\big(\frac{1}{4},\frac{1}{4},\frac{3}{4},\frac{3}{4},\frac{1}{6},\frac{5}{6}\big)$  &  5 & $(-5, 12, -18, 12, -5, 0)$&$B^2 A^3 B^{-3}$ & Yes \\ 
 \hline
27 & $\big(\frac{1}{4},\frac{3}{4},\frac{1}{5},\frac{2}{5},\frac{3}{5},\frac{4}{5}\big)$  & 7 & $(-7, 13, -22, 13, -7, 0)$& $A^4 B^4 A (AB)^{-1}$ &  Yes\\  
\hline
28 &  $\big(\frac{1}{4},\frac{3}{4},\frac{1}{6},\frac{5}{6},\frac{1}{6},\frac{5}{6}\big)$  & 4 & $(-4,11,-16,11,-4,0)$ & $BA^{-1}B^{6}A$ & Yes  \\ 
 \hline
29 &   $\big(\frac{1}{4},\frac{3}{4},\frac{1}{8},\frac{3}{8},\frac{5}{8},\frac{7}{8}\big)$ & 6 & $(-6, 14, -20, 14, -6, 0)$& $A^2 B^2 A^{-1} B^4 A B^{-1}A^{2}$ &  Yes\\  
\hline
30 &  $\big(\frac{1}{4},\frac{3}{4},\frac{1}{10},\frac{3}{10},\frac{7}{10},\frac{9}{10}\big) $ & 5 & $(-5, 13, -18, 13, -5, 0)$& $AB^2 A^{-1}B^{3}AB^{-1}$&Yes\\  
\hline
31 & $\big(\frac{1}{4},\frac{3}{4},\frac{1}{12},\frac{5}{12},\frac{7}{12},\frac{11}{12}\big)$ & 6 &$(-6, 15, -20, 15, -6, 0)$ & ?? &  ??   \\  
\hline
32 &  $\big(\frac{1}{5},\frac{2}{5},\frac{3}{5},\frac{4}{5},\frac{1}{6},\frac{5}{6}\big)$ & 6 & $(-6,14,-21,14,-6,0)$& $AB^{5}$ & Yes\\ 
 \hline
33 &   $\big(\frac{1}{6},\frac{5}{6},\frac{1}{6},\frac{5}{6},\frac{1}{6},\frac{5}{6}\big)$ & 3 & $(-3,9,-13,9,-3,0)$ &$A^2 B^4$ & Yes \\  
\hline
34 &  $\big(\frac{1}{6},\frac{5}{6},\frac{1}{8},\frac{3}{8},\frac{5}{8},\frac{7}{8}\big)$ & 5 &$(-5,14,-20,14,-5,0)$ &$B^{-4}AB^{4}$ & Yes\\
 \hline
35 &  $\big(\frac{1}{6},\frac{5}{6},\frac{1}{10},\frac{3}{10},\frac{7}{10},\frac{9}{10}\big)$ & 4 & $(-4,12,-17,12,-4,0)$& $B^{4} A$& Yes\\ 
 \hline
36 & $\big(\frac{1}{6},\frac{5}{6},\frac{1}{12},\frac{5}{12},\frac{7}{12},\frac{11}{12}\big)$ & 5 & $(-5,15,-21,15,-5,0)$& $BA^{-1}B^{6}AB^{-5}$ & Yes\\  
\hline
37 &  $\big(\frac{1}{7},\frac{2}{7},\frac{3}{7},\frac{4}{7},\frac{5}{7},\frac{6}{7}\big)$&7 & $(-7, 14, -21, 14, -7, 0)$& ?? & ?? \\ 
 \hline
38 & $\big(\frac{1}{9},\frac{2}{9},\frac{4}{9},\frac{5}{9},\frac{7}{9},\frac{8}{9}\big)$   & 6&  $(-6, 15, -21, 15, -6, 0)$&?? & ??\\  
\hline
39 &  $\big(\frac{1}{14},\frac{3}{14},\frac{5}{14},\frac{9}{14},\frac{11}{14},\frac{13}{14}\big)$  & 5& $(-5, 14, -19, 14, -5, 0)$& ??&?? \\ 
 \hline
40 &  $\big(\frac{1}{18},\frac{5}{18},\frac{7}{18},\frac{11}{18},\frac{13}{18},\frac{17}{18}\big)$ & 6 & $(-6, 15, -19, 15, -6, 0)$ & $A^4 B^4 A (A^{2}B)^{-1}$ & Yes\\ 
\hline
\hline
\end{longtable}
\end{center}
}

\begin{rmk}
In the $22$ cases whose arithmeticity or thinness is unknown, we have applied the SAGE program, written in Section \ref{se:sage}, with $k=15$: this shows in particular that, if one of the corresponding groups is arithmetic and a $\gamma$ satisfying the hypotheses of Proposition~\ref{Proposition} exists, then such a $\gamma$ has to be written as a product of at least $16$ matrices in $\{A,B,A^{-1},B^{-1}\}$. Among them, there are also $5$ cases (entries 1, 8, 15, 37, 38) for which arithmeticity cannot be proved through Proposition~\ref{Proposition}: in these cases, the gcd of the coordinates of the vector $v$ is larger than $2$, which implies that no $\gamma(v)$ can have $\pm 1,\pm 2$ in the last entry.
\end{rmk}
\section{Table B: More examples of arithmetic hypergeometric groups in $\Sp(6)$}\label{arithmeticexamples}
In this table we list all the $143$ pairs of the parameters $\alpha,\beta$ for which the leading coefficients of the difference polynomials $f-g$ have absolute values bigger than $2$, so the criterion of \cite[Theorem 1.1]{SV} cannot be applied in these cases but still the arithmeticity of the corresponding hypergeometric groups follows from Proposition \ref{Proposition}. Here the vector $v=(A^{-1}B-I)(e_6)$, $e_6=(0,0,0,0,0,1)$ and $\gamma\in \Gamma(f,g)$ for which the hypotheses of Proposition \ref{Proposition} are satisfied. Note that the values $lc(f-g)$, listed in the third column of Table A in the previous section, are nothing else but the first nonzero entry of the vectors $v$ and therefore we avoid to list these values in all the tables to follow, from Table B onwards.

{\begin{center}
\scriptsize\renewcommand{\arraystretch}{2}
\begin{longtable}{|c|c|c|c|c|c|}
\hline
 S.No. & $\alpha$ & $\beta$ & $v$  & $\gamma$  \\                                                                                                                                                        
\hline
\hline
1&$(0,0,0,0,\frac{1}{2},\frac{1}{2})$ & $\big(\frac{1}{3},\frac{1}{3},\frac{1}{3},\frac{2}{3},\frac{2}{3},\frac{2}{3}\big)$  &$(-5, -7, -3, -7, -5, 0)$ &$A B^{-1} A B^{2}A$  \\
\hline
 2&$(0,0,0,0,\frac{1}{2},\frac{1}{2})$ &$\big(\frac{1}{3},\frac{1}{3},\frac{2}{3},\frac{2}{3},\frac{1}{4},\frac{3}{4}\big)$   & $(-4, -5, 0, -5, -4, 0)$ &$ B^{4} A$ \\
 \hline
 3& $(0,0,0,0,\frac{1}{2},\frac{1}{2})$  & $\big(\frac{1}{3},\frac{2}{3},\frac{1}{4},\frac{3}{4},\frac{1}{4},\frac{3}{4}\big)$ &$(-3,-4,2,-4,-3,0)$ & $B^2$  \\
\hline
 4& $(0,0,0,0,\frac{1}{2},\frac{1}{2})$ & $\big(\frac{1}{3},\frac{2}{3},\frac{1}{5},\frac{2}{5},\frac{3}{5},\frac{4}{5}\big)$ &$(-4,-4,1,-4,-4,0)$ & $B^2A^{-1}B^3$  \\
\hline
5 & $(0,0,0,0,\frac{1}{2},\frac{1}{2})$ & $\big(\frac{1}{3},\frac{2}{3},\frac{1}{8},\frac{3}{8},\frac{5}{8},\frac{7}{8}\big)$ &$(-3,-2,4,-2,-3,0)$ & $B^2$  \\
\hline
6 & $(0,0,0,0,\frac{1}{2},\frac{1}{2})$ & $\big(\frac{1}{3},\frac{2}{3},\frac{1}{12},\frac{5}{12},\frac{7}{12},\frac{11}{12}\big)$ &$(-3,-1,5,-1,-3,0)$ & $B^2$ \\
\hline
7& $(0,0,0,0,\frac{1}{2},\frac{1}{2})$ & $\big(\frac{1}{4},\frac{3}{4},\frac{1}{5},\frac{2}{5},\frac{3}{5},\frac{4}{5}\big)$&$(-3, -3, 2, -3, -3, 0)$ & $A^3 B^{-2} A^{-1} B A^{-1}$ \\
\hline
\end{longtable}
\end{center}
}

\begin{rmk}
The arithmeticity of entry $142$, marked with an asterisk in the table above, has already been proved in~\cite{DFH} by computing its index. See entries 468 and 534 of \cite[Table 2]{DFH}, which represent the same group up to conjugation.
\end{rmk}

\section{Table-C: Examples of hypergeometric groups in $\Sp(6)$ for which arithmeticity follows from Singh and Venkataramana~\cite{SV}}
We list here the possible pairs of the parameters $\alpha,\beta$ for which the arithmeticity of the corresponding hypergeometric groups is determined by \cite[Theorem 1.1]{SV}. That is, for these cases  the leading coefficients of the difference polynomials $f-g$ have absolute values $1$ or $2$, so the arithmeticity of the corresponding hypergeometric groups follows directly from \cite[Theorem 1.1]{SV}.

\newpage
{\begin{center}
\scriptsize\renewcommand{\arraystretch}{2}
\begin{longtable}{|c|c|c|||c|c|c|}
\hline
 S.No. & $\alpha$ & $\beta$  & S.No.& $\alpha$ & $\beta$ \\                                                                                                                                                        
\hline
\hline
1 & $(0,0,0,0,\frac{1}{2},\frac{1}{2})$ &$\big(\frac{1}{3},\frac{2}{3},\frac{1}{4},\frac{3}{4},\frac{1}{6},\frac{5}{6}\big)$  & 2 & $(0,0,0,0,\frac{1}{2},\frac{1}{2})$ & $\big(\frac{1}{3},\frac{2}{3},\frac{1}{6},\frac{1}{6},\frac{5}{6},\frac{5}{6}\big)$  \\
\hline
 3 & $(0,0,0,0,\frac{1}{2},\frac{1}{2})$   & $\big(\frac{1}{3},\frac{2}{3},\frac{1}{10},\frac{3}{10},\frac{7}{10},\frac{9}{10}\big)$ & 4 & $(0,0,0,0,\frac{1}{2},\frac{1}{2})$ & $\big(\frac{1}{4},\frac{1}{4},\frac{1}{4},\frac{3}{4},\frac{3}{4},\frac{3}{4}\big)$  \\
 \hline
5 & $(0,0,0,0,\frac{1}{2},\frac{1}{2})$ & $\big(\frac{1}{4},\frac{1}{4},\frac{3}{4},\frac{3}{4},\frac{1}{6},\frac{5}{6}\big)$& 6 & $(0,0,0,0,\frac{1}{2},\frac{1}{2})$ & $\big(\frac{1}{4},\frac{3}{4},\frac{1}{8},\frac{3}{8},\frac{5}{8},\frac{7}{8}\big)$   \\
\hline
 7 & $(0,0,0,0,\frac{1}{2},\frac{1}{2})$   & $\big(\frac{1}{4},\frac{3}{4},\frac{1}{10},\frac{3}{10},\frac{7}{10},\frac{9}{10}\big) $ & 8 & $(0,0,0,0,\frac{1}{2},\frac{1}{2})$ &  $\big(\frac{1}{4},\frac{3}{4},\frac{1}{12},\frac{5}{12},\frac{7}{12},\frac{11}{12}\big)$  \\
 \hline
 9 & $(0,0,0,0,\frac{1}{2},\frac{1}{2})$ & $\big(\frac{1}{5},\frac{2}{5},\frac{3}{5},\frac{4}{5},\frac{1}{6},\frac{5}{6}\big)$  & 10 & $(0,0,0,0,\frac{1}{2},\frac{1}{2})$ & $\big(\frac{1}{6},\frac{5}{6},\frac{1}{6},\frac{5}{6},\frac{1}{6},\frac{5}{6}\big)$  \\
\hline
 11 & $(0,0,0,0,\frac{1}{2},\frac{1}{2})$ & $\big(\frac{1}{6},\frac{5}{6},\frac{1}{8},\frac{3}{8},\frac{5}{8},\frac{7}{8}\big)$ & 12 & $(0,0,0,0,\frac{1}{2},\frac{1}{2})$ & $\big(\frac{1}{6},\frac{5}{6},\frac{1}{12},\frac{5}{12},\frac{7}{12},\frac{11}{12}\big)$  \\
\hline
\hline
\end{longtable}
\end{center}
}

\section{Table-D: Other Open cases}
Following our study, we find that there are 86 pairs of the parameters $\alpha,\beta$ for which the leading coefficients of the difference polynomials $f-g$ have absolute values bigger than $2$,  and therefore the criterion of \cite[Theorem 1.1]{SV} cannot be applied in these cases. In addition,  we are also not able to find $\gamma\in\Gamma(f,g)$ which could satisfy the hypotheses of Proposition~\ref{Proposition}. Out of these, 22 are already listed in Table A in Section~\ref{sec:mum}. Here we list the remaining 64 pairs of parameters $\alpha,\beta$ for which the arithmeticity or thinness of the associated hypergeometric groups is still unknown.

All of them have been verified by the SAGE program written in Section~\ref{se:sage} for the values of $k$ up to $14$. If a $\gamma$ satisfying the hypotheses of Proposition~\ref{Proposition} exists for one of these cases, it must be a product of at least $15$ matrices in $\{A,B,A^{-1},B^{-1}\}$.

{\begin{center}
\scriptsize\renewcommand{\arraystretch}{2}
\begin{longtable}{|c|c|c|||c|c|c|}
\hline
 S.No. & $\alpha$ & $\beta$  & S.No.& $\alpha$ & $\beta$ \\                                                                                                                                                        
\hline
\hline
 1 & $\big(0,0,0,0,\frac{1}{2},\frac{1}{2}\big)$ & $\big(\frac{1}{3},\frac{1}{3},\frac{2}{3},\frac{2}{3},\frac{1}{6},\frac{5}{6}\big)$  & 2 & $(0,0,0,0,\frac{1}{3},\frac{2}{3})$ &  $\big(\frac{1}{2},\frac{1}{2},\frac{1}{2},\frac{1}{2},\frac{1}{4},\frac{3}{4}\big)$ \\
\hline
3 &$(0,0,0,0,\frac{1}{3},\frac{2}{3})$ & $\big(\frac{1}{2},\frac{1}{2},\frac{1}{2},\frac{1}{2},\frac{1}{6},\frac{5}{6}\big)$ & 4 &$(0,0,0,0,\frac{1}{3},\frac{2}{3})$ &$\big(\frac{1}{2},\frac{1}{2},\frac{1}{4},\frac{3}{4},\frac{1}{4},\frac{3}{4}\big)$ \\
\hline
 5 &$(0,0,0,0,\frac{1}{3},\frac{2}{3})$ & $\big(\frac{1}{2},\frac{1}{2},\frac{1}{5},\frac{2}{5},\frac{3}{5},\frac{4}{5}\big)$ &6 &$(0,0,0,0,\frac{1}{3},\frac{2}{3})$ & $\big(\frac{1}{2},\frac{1}{2},\frac{1}{8},\frac{3}{8},\frac{5}{8},\frac{7}{8}\big)$\\
 \hline
  7 & $(0,0,0,0,\frac{1}{3},\frac{2}{3})$&$\big(\frac{1}{2},\frac{1}{2},\frac{1}{10},\frac{3}{10},\frac{7}{10},\frac{9}{10}\big)$ & 8 &$(0,0,0,0,\frac{1}{3},\frac{2}{3})$ & $\big(\frac{1}{2},\frac{1}{2},\frac{1}{12},\frac{5}{12},\frac{7}{12},\frac{11}{12}\big)$   \\
\hline
9 & $(0,0,0,0,\frac{1}{3},\frac{2}{3})$& $\big(\frac{1}{7},\frac{2}{7},\frac{3}{7},\frac{4}{7},\frac{5}{7},\frac{6}{7}\big)$ & 10 &$(0,0,0,0,\frac{1}{3},\frac{2}{3})$ & $\big(\frac{1}{9},\frac{2}{9},\frac{4}{9},\frac{5}{9},\frac{7}{9},\frac{8}{9}\big)$\\
\hline
11 & $(0,0,0,0,\frac{1}{4},\frac{3}{4})$ & $\big(\frac{1}{2},\frac{1}{2},\frac{1}{2},\frac{1}{2},\frac{1}{3},\frac{2}{3}\big)$& 12 & $(0,0,0,0,\frac{1}{4},\frac{3}{4})$&$\big(\frac{1}{2},\frac{1}{2},\frac{1}{3},\frac{1}{3},\frac{2}{3},\frac{2}{3}\big)$ \\
\hline
13 &$(0,0,0,0,\frac{1}{4},\frac{3}{4})$ & $\big(\frac{1}{2},\frac{1}{2},\frac{1}{3},\frac{2}{3},\frac{1}{6},\frac{5}{6}\big)$& 14 &$(0,0,0,0,\frac{1}{4},\frac{3}{4})$ & $\big(\frac{1}{2},\frac{1}{2},\frac{1}{5},\frac{2}{5},\frac{3}{5},\frac{4}{5}\big)$ \\
\hline
15 & $(0,0,0,0,\frac{1}{4},\frac{3}{4})$& $\big(\frac{1}{2},\frac{1}{2},\frac{1}{6},\frac{5}{6},\frac{1}{6},\frac{5}{6}\big)$& 16 &$(0,0,0,0,\frac{1}{4},\frac{3}{4})$ & $\big(\frac{1}{2},\frac{1}{2},\frac{1}{8},\frac{3}{8},\frac{5}{8},\frac{7}{8}\big)$ \\
\hline
17 & $(0,0,0,0,\frac{1}{4},\frac{3}{4})$& $\big(\frac{1}{2},\frac{1}{2},\frac{1}{10},\frac{3}{10},\frac{7}{10},\frac{9}{10}\big)$& 18 &$(0,0,0,0,\frac{1}{4},\frac{3}{4})$& $\big(\frac{1}{2},\frac{1}{2},\frac{1}{12},\frac{5}{12},\frac{7}{12},\frac{11}{12}\big)$\\
\hline
 19 &$(0,0,0,0,\frac{1}{4},\frac{3}{4})$ & $\big(\frac{1}{7},\frac{2}{7},\frac{3}{7},\frac{4}{7},\frac{5}{7},\frac{6}{7}\big)$&  20 & $(0,0,0,0,\frac{1}{4},\frac{3}{4})$& $\big(\frac{1}{9},\frac{2}{9},\frac{4}{9},\frac{5}{9},\frac{7}{9},\frac{8}{9}\big)$ \\
 \hline
 21 & $(0,0,0,0,\frac{1}{6},\frac{5}{6})$  & $\big(\frac{1}{2},\frac{1}{2},\frac{1}{2},\frac{1}{2},\frac{1}{3},\frac{2}{3}\big)$& 22 &$(0,0,0,0,\frac{1}{6},\frac{5}{6})$ & $\big(\frac{1}{2},\frac{1}{2},\frac{1}{3},\frac{1}{3},\frac{2}{3},\frac{2}{3}\big)$\\
 \hline
  23 &$(0,0,0,0,\frac{1}{6},\frac{5}{6})$ & $\big(\frac{1}{2},\frac{1}{2},\frac{1}{3},\frac{2}{3},\frac{1}{4},\frac{3}{4}\big)$ &24 & $(0,0,0,0,\frac{1}{6},\frac{5}{6})$& $\big(\frac{1}{2},\frac{1}{2},\frac{1}{4},\frac{3}{4},\frac{1}{4},\frac{3}{4}\big)$\\
  \hline
   25 & $(0,0,0,0,\frac{1}{6},\frac{5}{6})$&  $\big(\frac{1}{2},\frac{1}{2},\frac{1}{5},\frac{2}{5},\frac{3}{5},\frac{4}{5}\big)$ & 26 & $(0,0,0,0,\frac{1}{6},\frac{5}{6})$& $\big(\frac{1}{2},\frac{1}{2},\frac{1}{8},\frac{3}{8},\frac{5}{8},\frac{7}{8}\big)$ \\
   \hline
    27 &$(0,0,0,0,\frac{1}{6},\frac{5}{6})$& $\big(\frac{1}{2},\frac{1}{2},\frac{1}{10},\frac{3}{10},\frac{7}{10},\frac{9}{10}\big)$  & 28 &$(0,0,0,0,\frac{1}{6},\frac{5}{6})$ & $\big(\frac{1}{2},\frac{1}{2},\frac{1}{12},\frac{5}{12},\frac{7}{12},\frac{11}{12}\big)$ \\
    \hline
 29 & $(0,0,0,0,\frac{1}{6},\frac{5}{6})$& $\big(\frac{1}{3},\frac{1}{3},\frac{1}{3},\frac{2}{3},\frac{2}{3},\frac{2}{3}\big)$ & 30 & $(0,0,0,0,\frac{1}{6},\frac{5}{6})$& $\big(\frac{1}{3},\frac{1}{3},\frac{2}{3},\frac{2}{3},\frac{1}{4},\frac{3}{4}\big)$\\
     \hline
 31 & $(0,0,0,0,\frac{1}{6},\frac{5}{6})$& $\big(\frac{1}{3},\frac{2}{3},\frac{1}{5},\frac{2}{5},\frac{3}{5},\frac{4}{5}\big)$ & 32 & $(0,0,0,0,\frac{1}{6},\frac{5}{6})$ & $\big(\frac{1}{4},\frac{3}{4},\frac{1}{12},\frac{5}{12},\frac{7}{12},\frac{11}{12}\big)$ \\
 \hline
 33 &$(0,0,0,0,\frac{1}{6},\frac{5}{6})$ &  $\big(\frac{1}{7},\frac{2}{7},\frac{3}{7},\frac{4}{7},\frac{5}{7},\frac{6}{7}\big)$&   34 & $(0,0,0,0,\frac{1}{6},\frac{5}{6})$& $\big(\frac{1}{9},\frac{2}{9},\frac{4}{9},\frac{5}{9},\frac{7}{9},\frac{8}{9}\big)$  \\
  \hline
 35 &$(0,0,\frac{1}{3},\frac{2}{3},\frac{1}{4},\frac{3}{4})$  & $(\frac{1}{2},\frac{1}{2},\frac{1}{5},\frac{2}{5},\frac{3}{5},\frac{4}{5})$ & 36 &$(0,0,\frac{1}{3},\frac{2}{3},\frac{1}{6},\frac{5}{6})$ & $(\frac{1}{2},\frac{1}{2},\frac{1}{4},\frac{1}{4},\frac{3}{4},\frac{3}{4})$ \\
 \hline
  37 &$(0,0,\frac{1}{3},\frac{2}{3},\frac{1}{6},\frac{5}{6})$  & $(\frac{1}{2},\frac{1}{2},\frac{1}{5},\frac{2}{5},\frac{3}{5},\frac{4}{5})$ & 38 & $(0,0,\frac{1}{3},\frac{2}{3},\frac{1}{6},\frac{5}{6})$ & $(\frac{1}{2},\frac{1}{2},\frac{1}{8},\frac{3}{8},\frac{5}{8},\frac{7}{8})$\\
  \hline
 
\hline
39 &$(0,0,\frac{1}{3},\frac{2}{3},\frac{1}{6},\frac{5}{6})$  &  $(\frac{1}{7},\frac{2}{7},\frac{3}{7},\frac{4}{7},\frac{5}{7},\frac{6}{7})$& 40 & $(0,0,\frac{1}{4},\frac{1}{4},\frac{3}{4},\frac{3}{4})$& $\big(\frac{1}{2},\frac{1}{2},\frac{1}{3},\frac{1}{3},\frac{2}{3},\frac{2}{3}\big)$ \\
\hline
41 &$(0,0,\frac{1}{4},\frac{1}{4},\frac{3}{4},\frac{3}{4})$ & $(\frac{1}{2},\frac{1}{2},\frac{1}{5},\frac{2}{5},\frac{3}{5},\frac{4}{5})$ & 42 & $(0,0,\frac{1}{4},\frac{1}{4},\frac{3}{4},\frac{3}{4})$& $(\frac{1}{3},\frac{2}{3},\frac{1}{12},\frac{5}{12},\frac{7}{12},\frac{11}{12})$  \\
 \hline
 43 & $(0,0,\frac{1}{4},\frac{3}{4},\frac{1}{6},\frac{5}{6})$&$\big(\frac{1}{2},\frac{1}{2},\frac{1}{3},\frac{2}{3},\frac{1}{3},\frac{2}{3}\big)$ & 44 & $(0,0,\frac{1}{4},\frac{3}{4},\frac{1}{6},\frac{5}{6})$&  $(\frac{1}{2},\frac{1}{2},\frac{1}{5},\frac{2}{5},\frac{3}{5},\frac{4}{5})$  \\
 \hline
 45 & $(0,0,\frac{1}{4},\frac{3}{4},\frac{1}{6},\frac{5}{6})$& $(\frac{1}{2},\frac{1}{2},\frac{1}{8},\frac{3}{8},\frac{5}{8},\frac{7}{8})$ & 46 &$(0,0,\frac{1}{4},\frac{3}{4},\frac{1}{6},\frac{5}{6})$ &  $\big(\frac{1}{3},\frac{1}{3},\frac{1}{3},\frac{2}{3},\frac{2}{3},\frac{2}{3}\big)$\\
 \hline
  47 &$(0,0,\frac{1}{4},\frac{3}{4},\frac{1}{6},\frac{5}{6})$ & $(\frac{1}{7},\frac{2}{7},\frac{3}{7},\frac{4}{7},\frac{5}{7},\frac{6}{7})$&  48 & $(0,0,\frac{1}{5},\frac{2}{5},\frac{3}{5},\frac{4}{5})$&  $\big(\frac{1}{2},\frac{1}{2},\frac{1}{3},\frac{1}{3},\frac{2}{3},\frac{2}{3}\big)$\\
  
 \hline
 49 & $(0,0,\frac{1}{6},\frac{1}{6},\frac{5}{6},\frac{5}{6})$& $\big(\frac{1}{2},\frac{1}{2},\frac{1}{3},\frac{1}{3},\frac{2}{3},\frac{2}{3}\big)$ & 50 &$(0,0,\frac{1}{6},\frac{1}{6},\frac{5}{6},\frac{5}{6})$ &  $(\frac{1}{2},\frac{1}{2},\frac{1}{5},\frac{2}{5},\frac{3}{5},\frac{4}{5})$ \\
 \hline
 51 &$(0,0,\frac{1}{6},\frac{1}{6},\frac{5}{6},\frac{5}{6})$ & $(\frac{1}{2},\frac{1}{2},\frac{1}{8},\frac{3}{8},\frac{5}{8},\frac{7}{8})$& 52 & $(0,0,\frac{1}{6},\frac{1}{6},\frac{5}{6},\frac{5}{6})$& $(\frac{1}{2},\frac{1}{2},\frac{1}{12},\frac{5}{12},\frac{7}{12},\frac{11}{12})$  \\
 \hline
 53 & $(0,0,\frac{1}{6},\frac{1}{6},\frac{5}{6},\frac{5}{6})$& $\big(\frac{1}{3},\frac{1}{3},\frac{1}{3},\frac{2}{3},\frac{2}{3},\frac{2}{3}\big)$& 54 & $(0,0,\frac{1}{6},\frac{1}{6},\frac{5}{6},\frac{5}{6})$&$\big(\frac{1}{3},\frac{1}{3},\frac{2}{3},\frac{2}{3},\frac{1}{4},\frac{3}{4}\big)$ \\
 \hline
 55 & $(0,0,\frac{1}{6},\frac{1}{6},\frac{5}{6},\frac{5}{6})$&$(\frac{1}{7},\frac{2}{7},\frac{3}{7},\frac{4}{7},\frac{5}{7},\frac{6}{7})$ &  56 &$(0,0,\frac{1}{6},\frac{1}{6},\frac{5}{6},\frac{5}{6})$ &$(\frac{1}{9},\frac{2}{9},\frac{4}{9},\frac{5}{9},\frac{7}{9},\frac{8}{9})$ \\
 \hline
  57 & $(0,0,\frac{1}{8},\frac{3}{8},\frac{5}{8},\frac{7}{8})$&$(\frac{1}{2},\frac{1}{2},\frac{1}{5},\frac{2}{5},\frac{3}{5},\frac{4}{5})$ & 58 & $(0,0,\frac{1}{8},\frac{3}{8},\frac{5}{8},\frac{7}{8})$ & $(\frac{1}{2},\frac{1}{2},\frac{1}{12},\frac{5}{12},\frac{7}{12},\frac{11}{12})$\\
  \hline
    59 &$(0,0,\frac{1}{10},\frac{3}{10},\frac{7}{10},\frac{9}{10})$  & $(\frac{1}{2},\frac{1}{2},\frac{1}{5},\frac{2}{5},\frac{3}{5},\frac{4}{5})$&  60 & $(0,0,\frac{1}{10},\frac{3}{10},\frac{7}{10},\frac{9}{10})$ & $(\frac{1}{2},\frac{1}{2},\frac{1}{12},\frac{5}{12},\frac{7}{12},\frac{11}{12})$\\
  \hline
   61 & $(0,0,\frac{1}{10},\frac{3}{10},\frac{7}{10},\frac{9}{10})$ & $(\frac{1}{7},\frac{2}{7},\frac{3}{7},\frac{4}{7},\frac{5}{7},\frac{6}{7})$ & 62 & $(0,0,\frac{1}{12},\frac{5}{12},\frac{7}{12},\frac{11}{12})$& $\big(\frac{1}{3},\frac{2}{3},\frac{1}{4},\frac{3}{4},\frac{1}{4},\frac{3}{4}\big)$\\
 \hline
  63 &  $(\frac{1}{3},\frac{1}{3},\frac{1}{3},\frac{2}{3},\frac{2}{3},\frac{2}{3})$&$\big(\frac{1}{6},\frac{1}{6},\frac{1}{6},\frac{5}{6},\frac{5}{6},\frac{5}{6}\big)$ & 64 & $(\frac{1}{3},\frac{1}{3},\frac{2}{3},\frac{2}{3},\frac{1}{3},\frac{2}{3})$ &$\big(\frac{1}{9},\frac{2}{9},\frac{4}{9},\frac{5}{9},\frac{7}{9},\frac{8}{9}\big)$ \\
 \hline
 \hline
\end{longtable}
\end{center}
}

\section*{Acknowledgements}
The first author would like to thank Albrecht Klemm for pointing out the possible connections of the various members of the family  listed in Table A with Calabi-Yau 5-folds and  referring to an important article~\cite{GMP}. 

The first and the third author take this opportunity to thank Wadim Zudilin for several discussions on the subject during their visit to MPIM, Bonn on various occasions. The first author extends his thanks to Peter Sarnak for interesting conversations about hypergeometric groups and  his encouragement at the initial stage of this work during his visit at CIRM, Luminy in December 2016.    

The work of the first and the second author is financially supported by ERC Consolidator grant 648329 (GRANT). The work of the third author is supported in part by the DST-INSPIRE Faculty Fellowship No.  DST/INSPIRE/04/2015/000794 and the SEED Grant No.  RD/0515-IRCCSH0-035 (IITBombay). 

\nocite{}
\bibliographystyle{abbrv}
\bibliography{BDSS}
\end{document}